\documentclass[twoside]{article} 
\usepackage{graphicx}
\date{} 
\title{The asymptotic expansion of a function due to L.L.~Karasheva}
\author{\sc R. B.\ Paris \\
{\em Division of Computing and Mathematics,} \\
{\em Abertay University, Dundee DD1 1HG, UK}}
\begin{document}
\def\f#1#2{\mbox{${\textstyle \frac{#1}{#2}}$}}
\def\dfrac#1#2{\displaystyle{\frac{#1}{#2}}}
\def\boldal{\mbox{\boldmath $\alpha$}}
\newcommand{\bee}{\begin{equation}}
\newcommand{\ee}{\end{equation}}
\newcommand{\sa}{\sigma}
\newcommand{\ka}{\kappa}
\newcommand{\al}{\alpha}
\newcommand{\la}{\lambda}
\newcommand{\ga}{\gamma}
\newcommand{\eps}{\epsilon}
\newcommand{\om}{\omega}
\newcommand{\fr}{\frac{1}{2}}
\newcommand{\fs}{\f{1}{2}}
\newcommand{\g}{\Gamma}
\newcommand{\br}{\biggr}
\newcommand{\bl}{\biggl}
\newcommand{\ra}{\rightarrow}
\newcommand{\gtwid}{\raisebox{-.8ex}{\mbox{$\stackrel{\textstyle >}{\sim}$}}}
\newcommand{\ltwid}{\raisebox{-.8ex}{\mbox{$\stackrel{\textstyle <}{\sim}$}}}
\renewcommand{\topfraction}{0.9}
\renewcommand{\bottomfraction}{0.9}
\renewcommand{\textfraction}{0.05}
\newcommand{\mcol}{\multicolumn}
\date{}
\maketitle
\pagestyle{myheadings}
\markboth{\hfill \sc R. B.\ Paris  \hfill}
{\hfill \sc Asymptotics of $F_{n,\sa}(x;\mu)$\hfill}
\begin{abstract}
We consider the asymptotic expansion for $x\to\pm\infty$ of the entire function 
\[F_{n,\sa}(x;\mu)=\sum_{k=0}^\infty \frac{\sin\,(n\gamma_k)}{\sin \gamma_k}\,\frac{x^k}{k! \Gamma(\mu-\sigma k)},\quad \gamma_k=\frac{(k+1)\pi}{2n}\]
for $\mu>0$, $0<\sigma<1$ and $n=1, 2, \ldots\ $. When $\sigma=\alpha/(2n)$, with $0<\alpha<1$, this function was recently introduced by L.L. Karasheva [{\it J. Math. Sciences}, {\bf 250} (2020) 753--759] as a solution of a fractional-order partial differential equation. 

By expressing $F_{n,\sa}(x;\mu)$ as a finite sum of Wright functions,
we employ the standard asymptotics of integral functions of hypergeometric type to determine its asymptotic expansion.
This is found to depend critically on the parameter $\sa$ (and to a lesser extent on the integer $n$). Numerical results are presented to illustrate the accuracy of the different expansions obtained.
\vspace{0.3cm}

\vspace{0.1cm}
 
\end{abstract}

\vspace{0.3cm}

\noindent $\,$\hrulefill $\,$

\vspace{0.3cm}
\begin{center}
{\bf 1.\ Introduction}
\end{center}
\setcounter{section}{1}
\setcounter{equation}{0}
\renewcommand{\theequation}{\arabic{section}.\arabic{equation}}
In a recent paper, L.L. Karasheva \cite{K} introduced the entire function 
\bee\label{e11}
\Theta_{n,\al}(x;\mu):=\sum_{k=0}^\infty \frac{\sin\,(n\gamma_k)}{\sin \gamma_k}\,\frac{x^k}{k! \g(\mu-\frac{\al k}{2n})},\qquad \gamma_k:=\frac{(k+1)\pi}{2n},
\ee
where $\mu>0$,  $0<\al<1$ and $n=1, 2, \ldots\ $ and throughout $x$ is a real variable. This function is of interest as it is involved in the fundamental solution of the differential equation
\[\frac{\partial^\al u}{\partial t^\al}+(-)^n \,\frac{\partial^{2n}u}{\partial x^{2n}}=f(x,t)\]
for positive integer $n$, where the derivative with respect to $t$ is the fractional derivative of order $\al$.
In the simplest case $n=1$, we have $\Theta_{1,\al}(x;\mu)=\phi(-\sa,\mu;x)$, $\sa:=\al/(2n)$, where $\phi(-\sa,\mu;x)$ is the Wright function 
\bee\label{e11b}
\phi(-\sa,\mu;x):=\sum_{k=0}^\infty \frac{x^k}{k! \g(\mu-\sa k)}\qquad (\sa<1),
\ee
which finds application as a fundamental solution of the diffusion-wave equation \cite{M}. 
Under the above assumptions on $n$ and $\al$ it follows that the parameter $\sa$ associated with (\ref{e11}) satisfies $0<\sa<\fs$. 

In this study, however, we shall allow the parameter $\sa$ to satisfy $0<\sa<1$ and consider the function
\bee\label{e11a}
F_{n,\sa}(x;\mu):=\sum_{k=0}^\infty \frac{\sin\,(n\gamma_k)}{\sin \gamma_k}\,\frac{x^k}{k! \g(\mu-\sa k)}\qquad(0<\sa<1),
\ee
which coincides with $\Theta_{n,\al}(x;\mu)$ when $\sa=\al/(2n)$.
From the well-known expansion
\[\frac{\sin\,(n\gamma_k)}{\sin \gamma_k}=\sum_{r=0}^{n-1} e^{i\gamma_k(2r-n+1)}=\sum_{r=0}^{n-1} e^{-i(k+1)\om_r}, \]
where
\bee\label{e12}
\om_r:=\frac{(n-2r-1)\pi}{2n}\qquad(0\leq r\leq n-1),
\ee
it follows that (\ref{e11a}) can be expressed as a finite sum of Wright functions defined in (\ref{e11b}) with rotated arguments (compare \cite[Eq.~(4)]{K})
\bee\label{e13}
F_{n,\sa}(x;\mu)=\sum_{r=0}^{n-1} e^{-i\om_r} \,\phi(-\sa,\mu;xe^{-i\om_r}).\hspace{0.7cm}
\ee 
We note that the extreme values of $\om_r$ satisfy $\om_{0}=-\om_{n-1}=(n-1)\pi/(2n)$, whence
$|\om_r|<\fs\pi$ for $0\leq r\leq n-1$.

The Wright function appearing in (\ref{e11b}) can be written alternatively as
\[\phi(-\sa,\mu;x)=\frac{1}{\pi}\sum_{k=0}^\infty \frac{x^k}{k!}\,\g(1-\mu+\sa k)\,\sin \pi(\mu-\sa k)\]
\[\hspace{4.2cm}=\frac{1}{2\pi}\bl\{e^{\pi i\vartheta} \Psi(xe^{\pi i\sa})+e^{-\pi i\vartheta} \Psi(xe^{-\pi i\sa})\br\}\hspace{2.2cm}\]
upon use of the reflection formula for the gamma function,
where $\vartheta:=\fs-\mu$ and the associated function $\Psi(z)$ is defined by
\bee\label{e15}
\hspace{2cm}\Psi(z):=\sum_{k=0}^\infty \frac{z^k}{k!}\,\g(\sa k+\delta)\qquad (0<\sa<1,\ \delta=1-\mu)
\ee
valid for $|z|<\infty$.

Hence we obtain the representation
\begin{eqnarray}
F_{n,\sa}(x;\mu)&=&\frac{1}{2\pi}\sum_{r=0}^{n-1} e^{-i\om_r} \Upsilon_r(\sa;x)\nonumber\\
&=&\frac{1}{\pi}\Re \bl\{\sum_{r=0}^{N-1}e^{-i\om_r} \Upsilon_r(\sa;x)+\Delta_n e^{\pi i\vartheta} \Psi(xe^{\pi i\sa})\br\},\label{e22}
\end{eqnarray}
where
\[\Upsilon_r(\sa;x):=e^{\pi i\vartheta} \Psi(xe^{\pi i\sa-i\om_r})+e^{-\pi i\vartheta} \Psi(xe^{-\pi i\sa-i\om_r})\]
and we have defined the quantities
\[N=\lfloor n/2\rfloor, \qquad \Delta_n=\left\{\begin{array}{ll}0 & (n\ \mbox{even})\\1 & (n\ \mbox{odd}).\end{array}\right.\]
The form (\ref{e22}) follows from the symmetry of the $\om_r$ in (\ref{e12}) (and the fact that $x$ is a real variable); 
the values of $\om_r$ for $0\leq r\leq N-1$ satisfy
\bee\label{e22a}
\{\om_0, \om_1,\ldots , \om_{N-1}\}=\bl\{\frac{(n-1)\pi}{2n}, \frac{(n-3)\pi}{2n}, \ldots , \frac{\pi}{2n}\, \epsilon_n\br\}, \qquad \epsilon_n=\left\{\begin{array}{ll} 1& (n\ \mbox{even})\\2& (n\ \mbox{odd}).\end{array}\right.
\ee

We shall use the representation in (\ref{e22}), with the above values of $\om_r$, to determine the asymptotic expansion of $F_{n,\sa}(x;\mu)$ for $x\to\pm\infty$ by application of the asymptotic theory of the integral function $\Psi(z)$. A summary of the expansion of $\Psi(z)$ for large $|z|$ is given in Section 2. The expansions of $F_{n,\sa}(x;\mu)$ for $x\to\pm\infty$ are given in Sections 3 and 4. A concluding section presents numerical results confirming the accuracy of the different expansions obtained. 

%which is based on the presentation discussed in \cite{P1, P2}.

\vspace{0.6cm}

\begin{center}
{\bf 2.\ The asymptotic expansion of $\Psi(z)$ for $|z|\to\infty$}
\end{center}
\setcounter{section}{2}
\setcounter{equation}{0}
\renewcommand{\theequation}{\arabic{section}.\arabic{equation}}
We first present the large-$|z|$ asymptotics of the function $\Psi(z)$ in (\ref{e15}) based on the presentation described in \cite[Section 4]{P1}; see also \cite[Section 4.2]{P2}, \cite[\S 2.3]{PK}. We introduce the following parameters:
\bee\label{a2}
\kappa=1-\sigma, \quad h=\sigma^\sigma , \quad \vartheta=\delta-\fs, \quad \delta=1-\mu,
\ee
together with the associated (formal) exponential and algebraic expansions
\bee\label{a3}
E(z):=Z^\vartheta e^Z\sum_{j=0}^\infty A_j(\sigma)Z^{-j},\qquad H(z):=\frac{1}{\sigma}\sum_{k=0}^\infty \frac{(-)^k}{k!} \g\bl(\frac{k+\delta}{\sigma}\br) z^{-(k+\delta)/\sigma},
\ee
where\footnote{The dependence of the coefficients $A_j(\sa)$ on the parameter $\delta$ is not indicated.} 
\bee\label{a4}
Z:=\ka(hz)^{1/\ka},\qquad A_0(\sigma)=(2\pi/\kappa)^{1/2}\,(\sigma/\kappa)^\vartheta.
\ee
Then, since $0<\kappa<1$, we obtain from \cite[p.~57]{PK} the large-$z$ expansion
\bee\label{a5}
\Psi(z)\sim \left\{\begin{array}{ll}E(z)+H(ze^{\mp\pi i}) & (|\arg\,z|\leq\fs\pi\kappa) \\
\\
H(ze^{\mp\pi i}) & (\fs\pi\kappa<|\arg\,z|\leq\pi),\end{array}\right.
\ee
where the upper or lower signs are chosen according as $\arg\,z>0$ or $\arg\,z<0$, respectively. 

The expansion $E(z)$ is exponentially large as $|z|\to\infty$ in the sector $|\arg\,z|<\fs\pi\kappa$, being oscillatory (multiplied by the algebraic factor $z^{\vartheta/\ka}$)  on the anti-Stokes lines $\arg\,z=\pm\fs\pi\ka$.
In the adjacent sectors $\fs\pi\ka<|\arg\,z|<\pi\ka$, the expansion $E(z)$ {\it continues to be present}, but is exponentially small reaching maximal subdominance relative to the algebraic expansion on the Stokes lines\footnote{On these rays $E(z)$ undergoes a Stokes phenomenon where it switches off in a smooth manner (see \cite[p.~67]{DLMF}).} $\arg\,z=\pm\pi\ka$.
In our treatment of $F_{n,\sa}(x;\mu)$ we shall not be concerned with exponentially small contributions, except in one special case when $x\to-\infty$ where the expansion of $F_{n,\sa}(x;\mu)$ is exponentially small. 

The first few normalised coefficients $c_j=A_j(\sa)/A_0(\sa)$ are \cite{P1, P2}:
\[c_0=1,\qquad c_1=\frac{1}{24\sigma}\{2+7\sigma+2\sigma^2-12\delta(1+\sigma)+12\delta^2\},\]
\[c_2=\frac{1}{1152\sigma^2}\{4+172\sigma+417\sigma^2+172\sigma^3+4\sigma^4-24\delta(6+41\sigma+41\sigma^2+6\sigma^3)\]
\[\hspace{6cm}+120\delta^2(4+11\sigma+4\sigma^2)-480\delta^3(1+\sigma)+144\delta^4\},\]
\[c_3=\frac{1}{414720 \sigma^3}\{(-1112 + 9636 \sigma + 163734 \sigma^2 + 336347 \sigma^3 + 
  163734 \sigma^4 + 9636 \sigma^5\]
  \[ - 1112 \sigma^6)-\delta(3600 + 220320 \sigma + 
  929700 \sigma^2 + 929700 \sigma^3  + 220320 \sigma^4  + 3600 \sigma^5)\] 
  \[+ \delta^2(65520  + 715680 \sigma + 1440180 \sigma^2 + 715680 \sigma^3 + 
  65520 \sigma^4)\] 
  \[ - \delta^3(161280  + 816480 \sigma  + 816480 \sigma^2  +
  161280 \sigma^3)\]
  \bee\label{ecj} +\delta^4 (151200  + 378000 \sigma + 151200 \sigma^2) - 
  60480\delta^5(1 + \sigma) + 8640 \delta^6\}.
  \ee

In addition to the Stokes lines $\arg\,z=\pm\pi\ka$, where $E(z)$ is maximally subdominant relative to the algebraic expansion, the positive real axis is also a Stokes line. Here the algebraic expansion is maximally subdominant relative to $E(z)$. As the positive real axis is crossed from the upper to the lower half plane the factor $e^{-\pi i}$ appearing in $H(ze^{-\pi i})$ changes to $e^{\pi i}$, and vice versa. The details of this transition will not be considered here; see \cite[Eq. (3.17)]{P92} for the case of the confluent hypergeometric function ${}_1F_1(a;b;z)$.
\vspace{0.6cm}

\begin{center}
{\bf 3.\ The asymptotic expansion of $F_{n,\sa}(x;\mu)$ for $x\to+\infty$}
\end{center}
\setcounter{section}{3}
\setcounter{equation}{0}
\renewcommand{\theequation}{\arabic{section}.\arabic{equation}}

\noindent{\bf 3.1}\ {\it Asymptotic character as a function of $\sa$}
\vspace{0.3cm}

Let us denote the arguments of the $\Psi$ functions appearing in (\ref{e22}) by
\[z_r^\pm=x \exp\,[i\phi_r^\pm],\qquad \phi_r^\pm=\pm\pi \sa-\om_r.\]
The representation of the asymptotic structure of the functions $\Psi(z_r^\pm)$ is illustrated in Fig.~1 for different values of $\sa$. The figures show the rays $\arg\,z=\pm\pi\sa$ and the anti-Stokes lines (dashed lines) $\arg\,z=\pm\fs\pi\ka$. In the case $\sa=\f{2}{3}$, the exponentially large sector is $|\arg\,z|<\f{1}{6}\pi$ and it is seen from Fig.~1(a) that the arguments $z_r^\pm$ for $0\leq r\leq N-1$ and $xe^{\pm\pi i\sa}$ all lie in the domain where $\Psi(z)$ has an algebraic expansion; this conclusion applies {\it a fortiori} when $\f{2}{3}<\sa<1$. When $\sa=\fs$, the exponentially large sector is $|\arg\,z|<\f{1}{4}\pi$; when $n=2$ we have $\om_0=\f{1}{4}\pi$ so that $z_0^+$ is situated on the boundary of the exponentially large sector. Other values of $n\geq 3$ will have some $z_r^+$ inside this sector, whereas the $z_r^-$ are in the algebraic sector for $n\geq 2$. Similarly, the case $\sa=\f{1}{3}$, where the rays $\arg\,z=\pm\pi\sa$ and $\arg\,z=\pm\fs\pi\ka$ coincide,  has all the $z_r^+$ situated in the exponentially large sector, with the $z_r^-$ situated in the algebraic domain. Finally, when $\sa=\f{1}{6}$ the exponentially large sector $|\arg\,z|<\f{5}{12}\pi$ encloses the rays $\arg\,z=\pm\pi\sa$, with the result that all the $z_r^+$ lie in the exponentially large sector, whereas the $z_r^-$ lie in the algebraic domain (except when $n=2$ when $z_0^-$ lies on the lower boundary of the exponentially large sector).

\begin{figure}[th]
	\begin{center}	{\tiny($a$)}\includegraphics[width=0.26\textwidth]{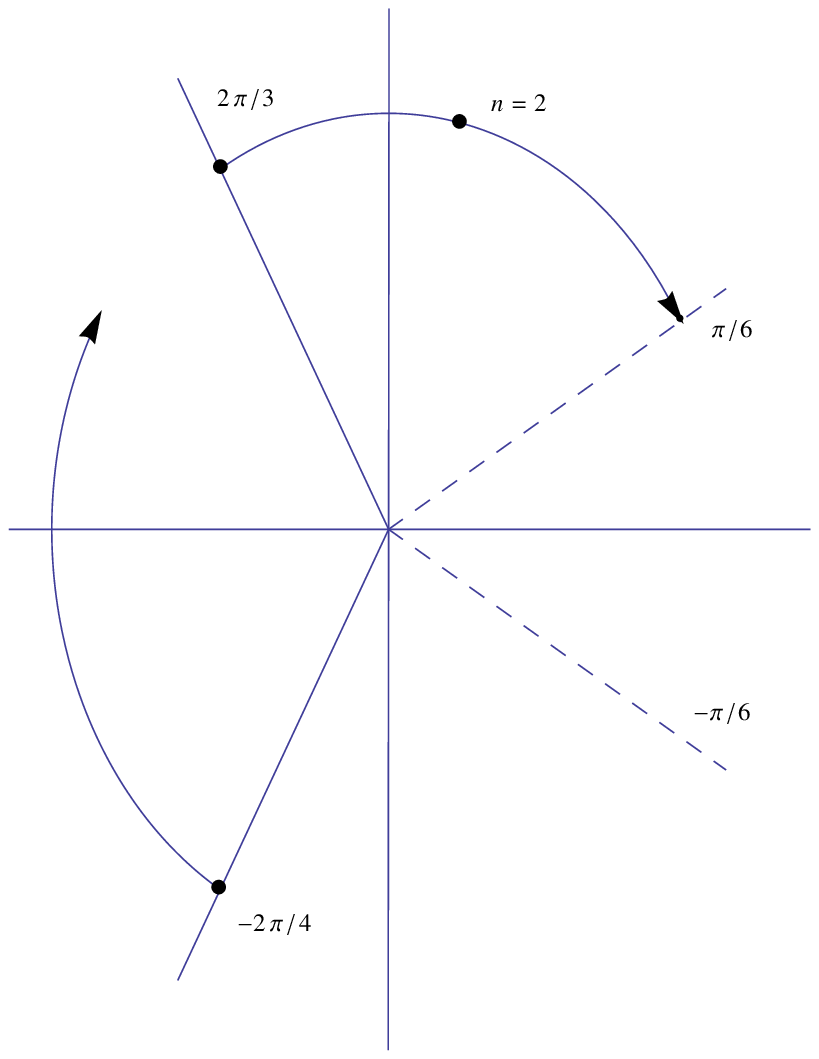}
	\qquad
	{\tiny($b$)}\includegraphics[width=0.3\textwidth]{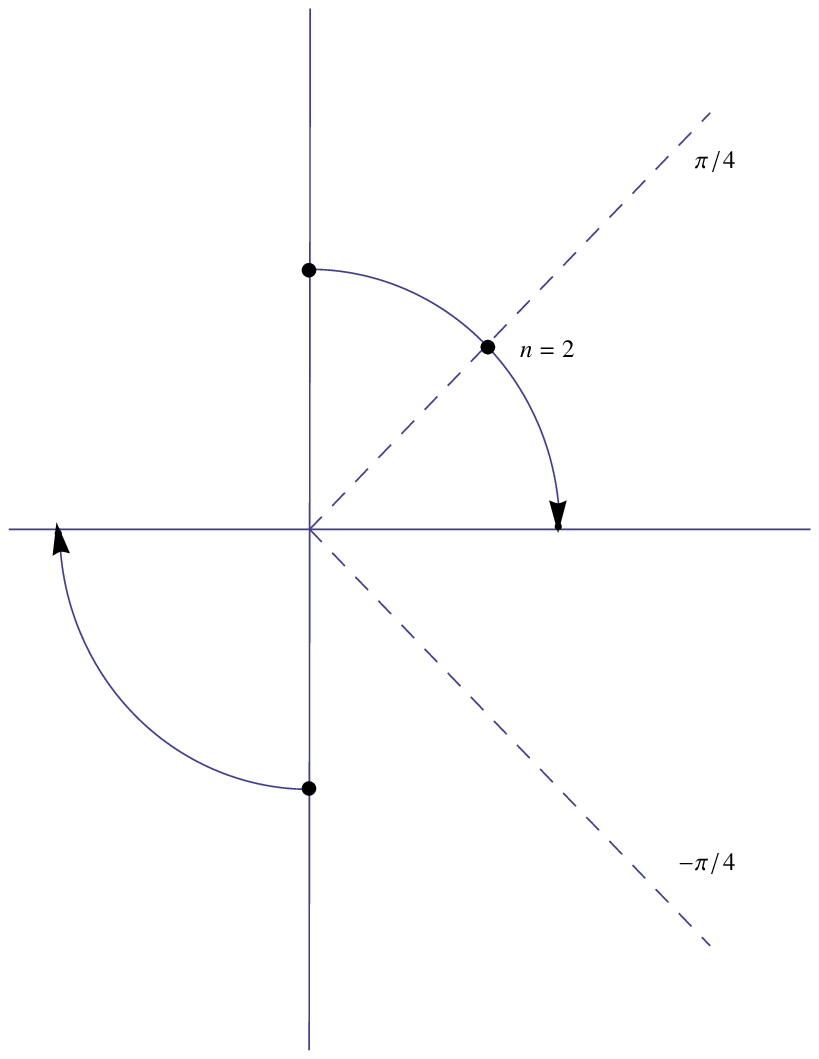} 
	\bigskip
	
		{\tiny($c$)}\includegraphics[width=0.3\textwidth]{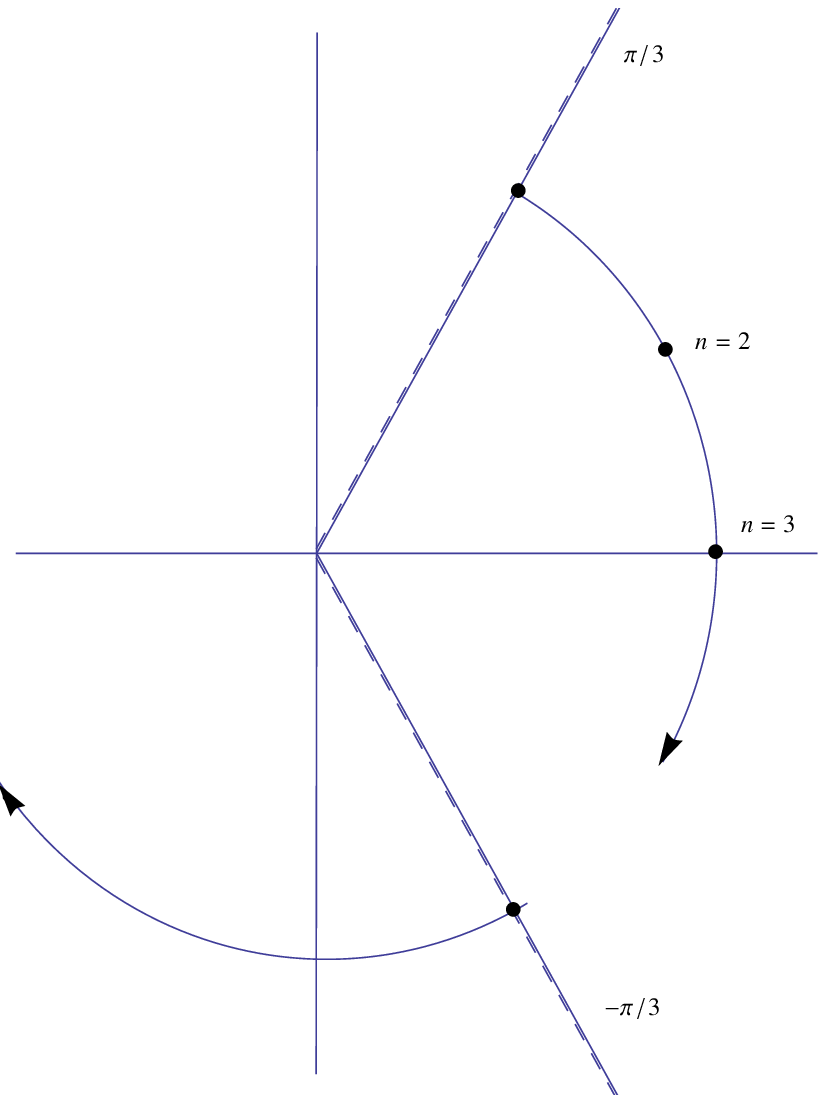}
	\qquad
	{\tiny($d$)}\includegraphics[width=0.3\textwidth]{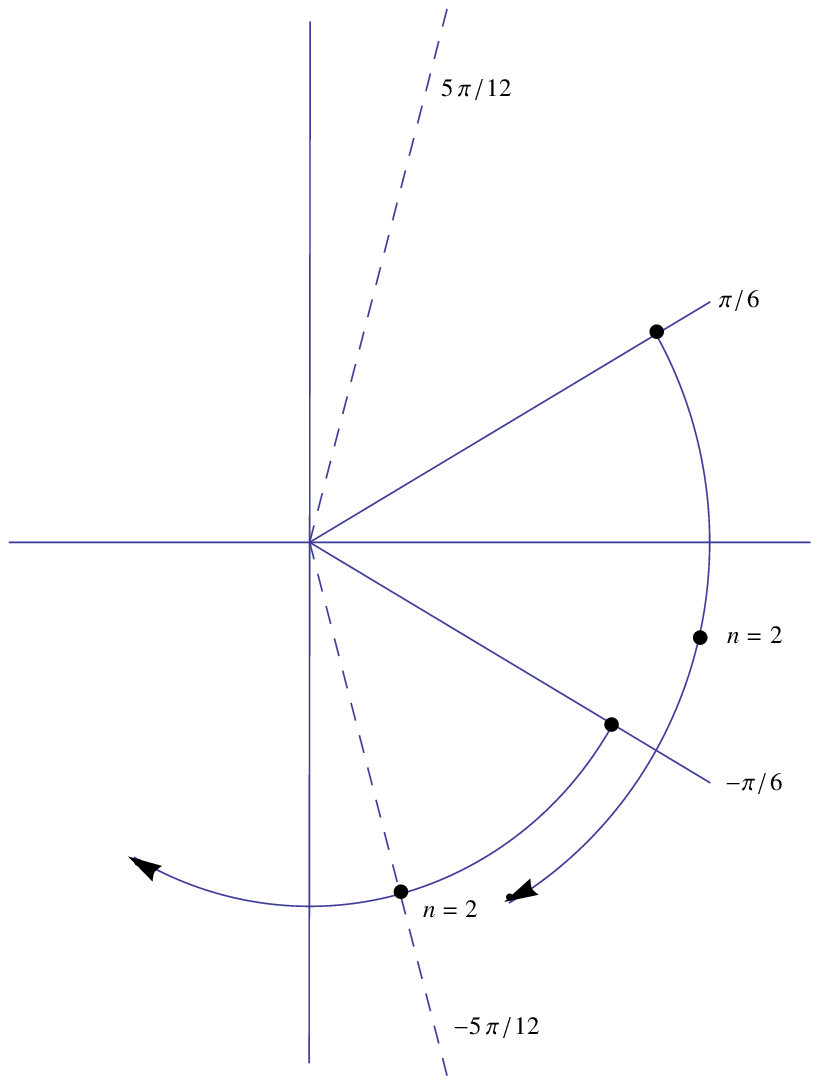} \\
	
\caption{\small{Diagrams representing the rays $\arg\,z=\pm\pi\sa$ and the boundaries of the exponentially large sector (shown by dashed rays) $|\arg\,z|<\fs\pi\ka$, $\ka=1-\sa$ for (a) $\sa=2/3$, (b) $\sa=1/2$, (c) $\sa=1/3$ and (d) $\sa=1/6$. Outside the exponentially large sector the expansion of $\Psi(z)$ is algebraic in character. The circular quadrants represent the range of the arguments $\arg\,z=\pm\pi\sa-\om_r$ for $0\leq r\leq \lfloor n/2\rfloor-1$, with $n\geq 2$ and the arrow-head corresponds to $n=\infty$. When $\sa=1/3$ the rays $\arg\,z=\pm\pi\sa$ and $\arg\,z=\pm\fs\pi\ka$ coincide.}}
\end{center}
\end{figure}

To summarise, we have the following asymptotic character of $F_{n,\sa}(x;\mu)$ when $x\to+\infty$ as a function of the parameter $\sa$:
\bee\label{e24}
\left.\begin{array}{ll}
0<\sa<\fs & \mbox{Exp.\ large}+ \mbox{Algebraic\ (for}\ n\geq 2)\\
\\
\fs\leq\sa<\f{2}{3}& \mbox{Exp,\ large}\ (\mbox{dependent\ on}\ n)+ \mbox{Algebraic}\  \\
\\
\f{2}{3}\leq\sa<1 & \mbox{Algebraic} \  (\mbox{for}\ n\geq 2). \end{array}\right\}
\ee
%where EL and A denote `exponentially large' and `algebraic' behaviour, respectively.
\vspace{0.3cm}

\noindent{\bf 3.2}\ {\it Asymptotic expansion}
\vspace{0.3cm}

\noindent
From (\ref{e22}) and (\ref{a3}), we have the algebraic expansion associated with $F_{n,\sa}(x;\mu)$ given by
\bee\label{e25}
{\bf H}(x)=\frac{1}{\sa}\sum_{k=0}^\infty \frac{x^{-K}}{k! \g(1-K)}\,\theta_{n,k}~,\qquad K:=\frac{k+\delta}{\sa},
\ee
where, with appropriate choices of the factors $e^{\pm\pi i}$ in $H(z)$,
\begin{eqnarray}
\theta_{n,k}&=&\frac{(-)^k}{\sin \pi K} \Re \bl\{\sum_{r=0}^{N-1} e^{\pi i\vartheta-i\om_r} (e^{\pi i\sa-i\om_r}\cdot e^{-\pi i})^{-K}+e^{-\pi i\vartheta-i\om_r}(e^{-\pi i\sa-i\om_r}\cdot e^{\pi i})^{-K}\nonumber\\
&&\hspace{8cm}+\Delta_n e^{\pi i\vartheta} (e^{\pi i\sa}\cdot e^{-\pi i})^{-K}\br\}\nonumber\\
&=&\frac{(-)^k}{\sin \pi K} \Re \bl\{\sum_{r=0}^{N-1} e^{(K-1)i\om_r}(e^{\pi i(\vartheta+\ka K)}+e^{-\pi i(\vartheta+\ka K)})+\Delta_n e^{\pi i(\vartheta+\ka K)}\br\}\nonumber\\
&=&\Re \bl\{2\sum_{r=0}^{N-1} e^{(K-1)i\om_r}+\Delta_n\br\},\label{e25a}
\end{eqnarray}
since $\cos \pi(\vartheta+\ka K)=\cos \pi(K-k-\fs)=(-)^k \sin \pi K$.

For the exponential component we introduce the quantities
\bee\label{e26a}
X=\ka(hx)^{1/\ka},\qquad \Phi_r^\pm=\pm\frac{\pi\vartheta}{\ka}-\om_r\bl(1+\frac{\vartheta}{\ka}\br)
\ee
and the formal asymptotic sum
\bee\label{e26b}
S(Xe^{i\Omega}):=\sum_{j=0}^\infty A_j(\sa) (Xe^{i\Omega/\ka})^{-j}.
\ee
Then, from (\ref{e22}) and (\ref{a3}), we have the exponential expansion in the form
\[{\bf E}(x)=\frac{X^\vartheta}{\pi}\Re \bl\{\sum_{r=0}^{N-1} \bl(\exp\,[Xe^{i\phi_r^+/\ka}+i\Phi_r^+]\, S(Xe^{i\phi_r^+})+\exp\,[Xe^{i\phi_r^-/\ka}+i\Phi_r^-]\,S(Xe^{i\phi_r^-})\br)\]
\bee\label{e26}
+\Delta_n \exp\,[Xe^{\pi i\sa/\ka}+\pi i\vartheta/\ka]\,S(Xe^{\pi i\sa})\br\}.
\ee

It is important to stress that only the exponential terms with $|\phi_r^\pm|\leq\fs\pi\ka$, that is those with
\[|\pm\pi\sa-\om_r|\leq\fs\pi\ka,\]
are to be retained in ${\bf E}(x)$ in (\ref{e26}). In addition, it is seen by inspection of Fig.~1 that the second term involving $S(Xe^{i\phi_r^-})$ does not contribute to ${\bf E}(x)$ when $\f{1}{3}\leq\sa<1$, since for this range of $\sa$ the ray $\arg\,z=-\pi\sa$ lies outside (or, when $\sa=\f{1}{3}$, on the lower boundary of) the exponentially large sector $|\arg\,z|<\fs\pi\ka$. Thus, when $\fs\leq\sa<\f{2}{3}$, the exponential expansion is significant if $\pi\sa-\om_0\leq\fs\pi\ka$; that is, if $n\geq n_0=1/(2-3\sa)$.

In summary, we have the following theorem:
\newtheorem{theorem}{Theorem}
\begin{theorem}$\!\!\!.$ \ The following expansion holds for $x\to+\infty$:
\[F_{n,\sa}(x;\mu)\sim\left\{\begin{array}{ll} {\bf E}(x)+{\bf H}(x) & (0<\sa<\fs;\ n\geq 2)\\
\\
{\bf E}(x)+{\bf H}(x) & (\fs\leq\sa<\f{2}{3};\ n\geq n_0)\\
\\
{\bf H}(x) & (\fs\leq\sa<\f{2}{3};\ n<n_0)\\
\\
{\bf H}(x) & (\f{2}{3}\leq \sa<1;\ n\geq 2),\end{array}\right.\]
where $n_0=1/(2-3\sa)$ and the exponential and algebraic expansions ${\bf E}(x)$ and ${\bf H}(x)$ are defined in (\ref{e25}) and (\ref{e26}).
\end{theorem}
\vspace{0.3cm}

\noindent{\bf 3.3}\ {\it Karasheva's estimate for $|\Theta_{n,\al}(x;\mu)|$}
\vspace{0.3cm}

\noindent
When $\sa=\al/(2n)<\fs$, we see from Theorem 1 that the dominant exponential expansion as $x\to+\infty$ corresponds to $r=0$ yielding
\[\Theta_{n,\al}(x;\mu)\sim \frac{A_0(\sa) X^\vartheta}{\pi} \Re \exp\,[X e^{i(\pi\sa-\om_0)/\ka+i\Phi_0^+}]\]
\[=\frac{A_0(\sa) X^\vartheta}{\pi} \exp\,[X\cos (\pi\sa-\om_0)/\ka)] \,\cos [X \sin (\pi\sa-\om_0)/\ka)+\Phi_0^+],\]
where
\[\frac{\pi\sa-\om_0}{\ka}=\frac{2n\pi\sa-(n-1)\pi}{2n-\al}=\frac{(\al+1-n)\pi}{2n-\al}~.\]

Thus we have the leading order estimate
\bee\label{e27}
\Theta_{n,\al}(x;\mu)\sim \frac{A_0(\sa) X^\vartheta}{\pi} \exp\,\bl[X\cos \bl(\frac{(n\!-\!1\!-\!\al)\pi}{2n\!-\!\al}\br)\br]
\cos\,\bl[X \sin \bl(\frac{(n\!-\!1\!-\!\al)\pi}{2n\!-\!\al}\br)-\Phi_0^+\br]
\ee 
as $x\to+\infty$. When expressed in our notation, Karasheva's estimate for $|\Theta_{n,\al}(x;\mu)|$ in \cite[\S 8]{K} agrees with (\ref{e27}) (when the second cosine term is replaced by 1), except that she did not give the value of the multiplicative constant $A_0(\sa)/\pi$ given in (\ref{a4}).
However, the presentation of her result as an upper bound is not evident due to the presence of possibly less dominant exponential expansions and also the subdominant algebraic expansion.

\vspace{0.6cm}

\begin{center}
{\bf 4.\ The expansion of $F_{n,\sa}(x;\mu)$ for $x\to-\infty$}
\end{center}
\setcounter{section}{4}
\setcounter{equation}{0}
\renewcommand{\theequation}{\arabic{section}.\arabic{equation}}
To examine the case of negative $x$ we replace $x$ by $e^{\mp\pi i}x$, with $x>0$, and use the fact that $\Psi(ze^{2\pi i})=\Psi(z)$ to find from (\ref{e22}) that
\bee\label{e31}
F_{n,\sa}(-x;\mu)=\frac{1}{\pi}\Re \bl\{\sum_{r=0}^{N-1}e^{-i\om_r} \Upsilon_r(-\ka;x)
%\bl(e^{\pi i\vartheta}\Psi(xe^{-\pi i\ka-i\om_r})+
%e^{-\pi i\vartheta}\Psi(xe^{\pi i\ka-i\om_r})\br)
+\Delta_n\,e^{\pi i\vartheta} \Psi(xe^{-\pi i\ka})\br\}.
\ee
The rays $\arg\,z=\pm\pi\sa$ in Fig.~1 are now replaced by the Stokes lines $\arg\,z=\pm\pi\ka$. The Stokes and anti-Stokes lines $\arg\,z=\pm\fs\pi\ka$ are illustrated in Fig.~2 when $0<\sa<\fs$ and $\fs<\sa<1$. In the sectors $\fs\pi\ka<|\arg\,z|<\pi\ka$, we recall that the exponential expansion $E(z)$ is still present but is exponentially small as $|z|\to\infty$.
\begin{figure}[th]
	\begin{center}	{\tiny($a$)}\includegraphics[width=0.3\textwidth]{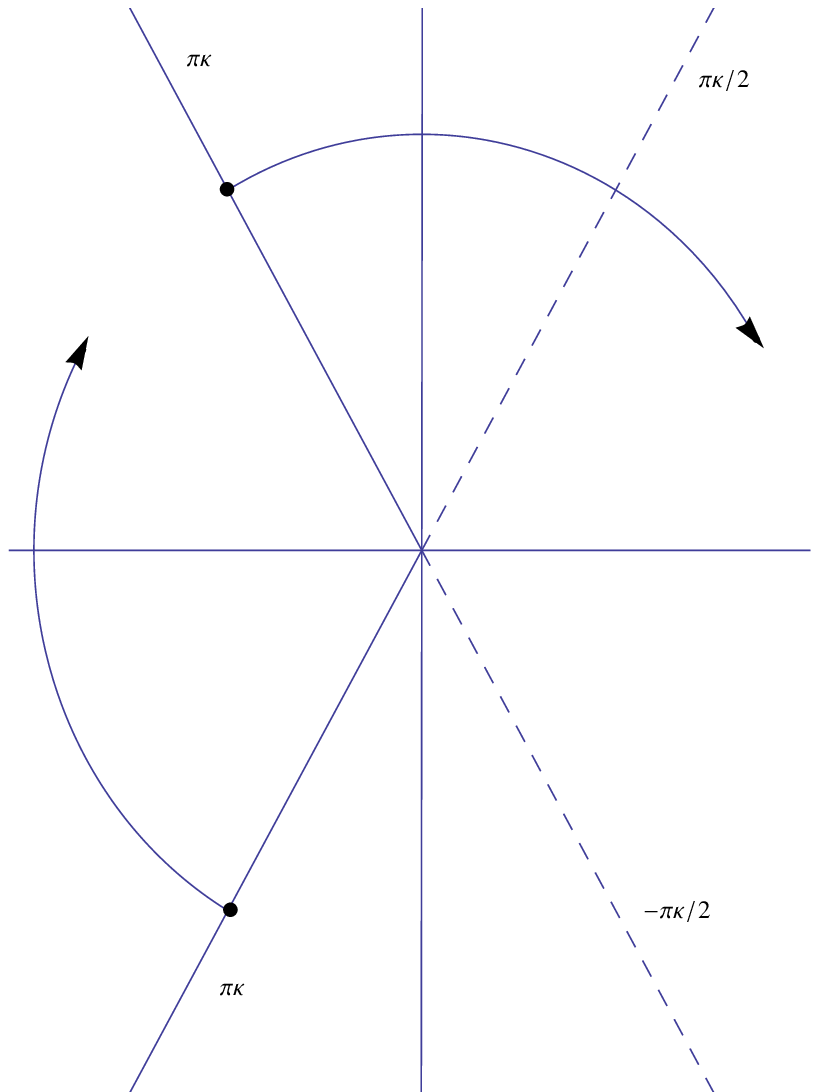}
	\qquad
	{\tiny($b$)}\includegraphics[width=0.3\textwidth]{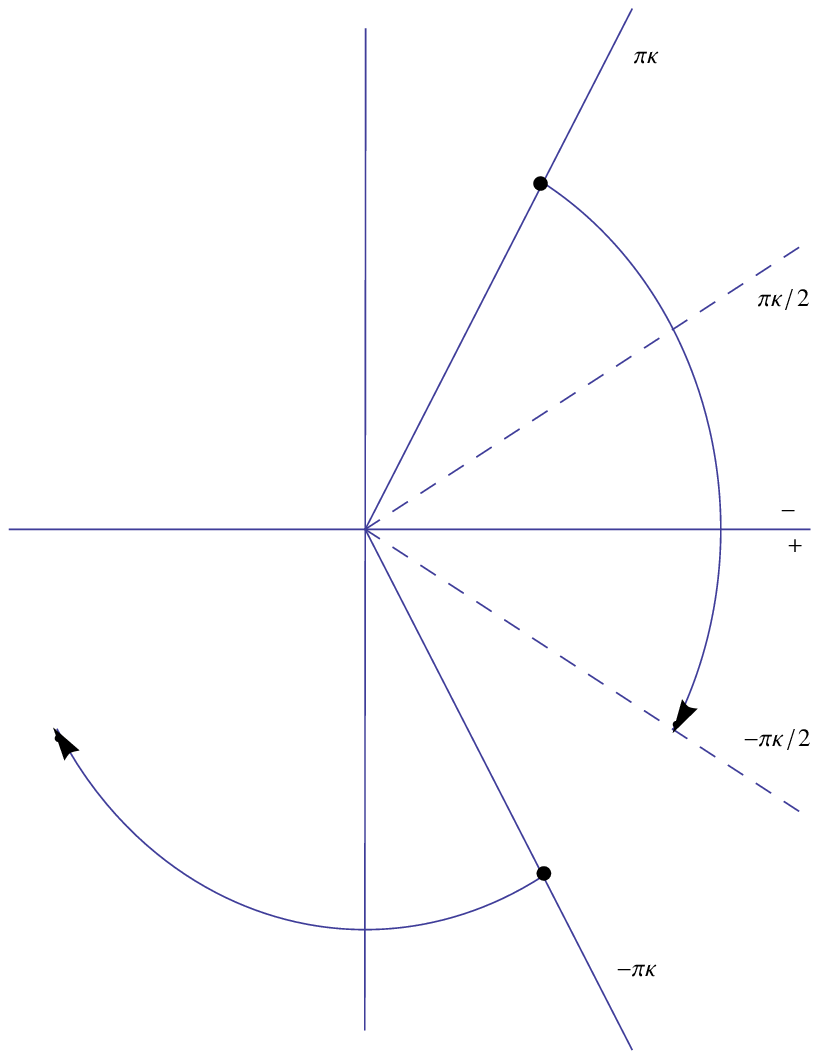} \\
	
\caption{\small{Diagrams representing the rays $\arg\,z=\pm\pi\ka$ and the boundaries of the exponentially large sector (shown by dashed rays) $|\arg\,z|<\fs\pi\ka$, $\ka=1-\sa$ for (a) $0<\sa<\fs$ and (b) $\fs<\sa<1$. The circular quadrants represent the range of the arguments $\arg\,z=\pm\pi\ka-\om_r$ for $0\leq r\leq N-1$ with the arrow-head corresponding to $n=\infty$. The $\pm$ signs in (b) denote the signs to be chosen in $H(z)$ on either side of the Stokes line $\arg\,z=0$.}}
\end{center}
\end{figure}

For the algebraic component of the expansion two cases arise when the argument $\pi\ka-\om_r$ of the second $\Psi$ function in $\Upsilon_r(-\ka;x)$  is either (i) positive or (ii) negative. In case (i) the algebraic expansion $H(z)$ does not encounter a Stokes phenomenon as its argument does not cross $\arg\,z=0$, whereas in case (ii) a Stokes phenomenon arises for those values of $r$ that make $\pi\ka-\om_r<0$.
In case (i), the algebraic component contains the factor inside the sum over $r$ in (\ref{e31})
\[e^{\pi i\vartheta}(e^{-\pi i\ka-i\om_r}\cdot e^{\pi i})^{-K}+e^{-\pi i\vartheta}(e^{\pi i\ka-i\om_r}\cdot e^{-\pi i})^{-K}\]
\[=e^{i\om_rK}(e^{\pi i(\vartheta-\sa K)}+e^{-\pi i(\vartheta-\sa K)})=2e^{i\om_r K} \cos \pi(k+\fs)\equiv 0\]
upon recalling the definition of $K$ in (\ref{e25}) and noting that $\delta-\vartheta=\fs$. Similarly, the final term involves the factor $\Re e^{\pi i\vartheta}(e^{-\pi i\ka}\cdot e^{\pi i})^{-K}=\cos \pi(\vartheta-\sa K)=0$. Thus the algebraic contribution to $F_{n,\sa}(-x;\mu)$ vanishes in case (i).

For case (ii) to apply, we require that $\pi\ka-\om_0<0$; that is, $n>n^*=1/(2\sa-1)$. Suppose that $\pi\ka-\om_r<0$ for $0\leq r\leq r_0$. Then the algebraic component resulting from the terms with $r\leq r_0$ becomes
\[\frac{1}{\pi\sa} \Re \bl\{\sum_{k=0}^\infty \frac{(-)^k \g(K)}{k!}\,x^{-K} \sum_{r=0}^{r_0}e^{(K-1)i\om_r}\bl(e^{\pi i\vartheta}(e^{-\pi i\ka}\cdot e^{\pi i})^{-K}+e^{-\pi i\vartheta}(e^{\pi i\ka}\cdot e^{\pi i})^{-K}\br)\br\}\]
\[=\frac{2}{\pi\sa}\Re \bl\{\sum_{k=0}^\infty \frac{(-)^k \g(K)}{k!}\,x^{-K} \sum_{r=0}^{r_0} e^{(K-1)i\om_r-\pi iK} \cos \pi(\vartheta-\sa K+\pi K)\br\},\]
where in the second term in round braces we have taken account of the Stokes phenomenon (the first term and that multiplied by $\Delta_n$ are unaffected).
Some routine algebra then produces the algebraic contribution
\bee\label{e32}
{\bf{\hat H}}(x):=\frac{2}{\sa}\sum_{k=0}^\infty \frac{x^{-K}}{k! \g(1-K)}\,{\hat \theta}_{n,k},\qquad {\hat \theta}_{n,k}:= \sum_{r=0}^{r_0} \cos\,\bl\{\pi K-(K-1)\om_r\br\}
\ee
when $n>n^*$, and\footnote{We avoid here consideration of the algebraic contribution when $\pi\ka-\om_r=0$, that is, on the Stokes line $\arg\,z=0$.} %In any case, the algebraic contribution is subdominant since ${\bf{\hat E}}(x)$ is exponentially large.} 
${\bf{\hat H}}(x)\equiv 0$ when $n<n^*$.

Reference to Fig.~2 shows that there is no exponential contribution to $F_{n,\sa}(-x;\mu)$ from the terms $\Psi(xe^{-\pi i\ka})$ and $\Psi(xe^{-\pi i\ka-i\om_r})$. From (\ref{a3}) and (\ref{e31}), we find the exponential expansion results from the terms $\Psi(xe^{\pi i\ka-i\om_r})$ and is given by
\bee\label{e33}
{\bf{\hat E}}(x):=\frac{X^\vartheta}{\pi} \Re \sum_{r=0}^{N-1} \exp\,[-Xe^{-i\om_r/\ka}-i\Phi] \ S(-Xe^{-i\om_r/\ka}),
\ee
where $X$ and the asymptotic sum $S$ are defined in (\ref{e26a}) and (\ref{e26b}) with $\Phi:=\om_r (1+\vartheta/\ka)$. For $\sa<\fs$ (when the algebraic expansion vanishes) the expansion of $F_{n,\sa}(-x;\mu)$ will be exponentially small as $|z|\to\infty$ provided $\pi\ka-\om_0>\fs\pi\ka$; that is, when $n<1/\sa$. If $n=1/\sa$, there is an exponentially oscillatory contribution and when $n>1/\sa$, the expansion is exponentially large.

To summarise we have the theorem:
\begin{theorem}$\!\!\!.$\ The following expansion holds for $x\to-\infty$:
\bee\label{e34}
F_{n,\sa}(-x;\mu)\sim\left\{\begin{array}{ll} {\bf{\hat E}}(x) & (0<\sa\leq\fs) \\
\\
{\bf{\hat E}}(x)+{\bf{\hat H}}(x) & (\fs<\sa<1),\end{array}\right.
\ee
where the exponential expansion ${\bf{\hat E}}(x)$ is defined in (\ref{e33}). This last expansion is exponentially small as $x\to-\infty$ when $0<\sa<\fs$ and $n<1/\sa$.
The algebraic expansion ${\bf{\hat H}}(x)$ is given by
%\[{\bf{\hat H}}(x):=\left\{\begin{array}{ll} \dfrac{2}{\sa}\sum_{k=0}^\infty \dfrac{x^{-K}}{k! \g(1-K)}\,{\hat \theta}_{n,k} & (n>n^*) \\
%\\
%0 & (n<n^*), \end{array}\right.\]
\[{\bf{\hat H}}(x):=\frac{2}{\sa}\sum_{k=0}^\infty \dfrac{x^{-K}}{k! \g(1-K)}\,{\hat \theta}_{n,k} \ \  (n>n^*),\qquad 0\ \ (n<n^*),\] 
where $n^*=1/(2\sa-1)$ and $K$, ${\hat\theta}_{n,k}$ are specified in (\ref{e25}) and (\ref{e32}).
\end{theorem}

\vspace{0.6cm}

\begin{center}
{\bf 5.\ Numerical results}
\end{center}
\setcounter{section}{5}
\setcounter{equation}{0}
\renewcommand{\theequation}{\arabic{section}.\arabic{equation}}
In this section we describe numerical calculations that support the expansions given in Theorems 1 and 2. The function $F_{n,\sa}(x;\mu)$ was evaluated using the expression in terms of Wright functions (valid for real $x$)
\bee\label{e41}
F_{n,\sa}(x;\mu)=2\Re \sum_{r=0}^{N-1} e^{i\om_r} \phi(-\sa,\mu; xe^{i\om_r})+\Delta_n \phi(-\sa,\mu;x),\qquad N=\lfloor n/2\rfloor,
\ee
which follows from (\ref{e13}) and the symmetry of $\om_r$. 
\begin{table}[h]
\caption{\footnotesize{The values of the exponential and algebraic expansions compared with $F_{n,\sa}(x;\mu)$ for large $x>0$ for different values of $\sa$ and $n$ when $\mu=3/4$ and $x=8$.}}
\begin{center}
\begin{tabular}{|c r|l|l|l|}
\hline
&&&&\\[-0.3cm]
%\mcol{1}{|c|}{} & \mcol{1}{c|}{$\gamma=0$} &\mcol{1}{c|}{$\gamma=2$}& \mcol{1}{c|}{$\gamma=-2$} \\
\mcol{1}{|c}{$\sa$} & \mcol{1}{c|}{} & \mcol{1}{c|}{$n=2$} &\mcol{1}{c|}{$n=3$}& \mcol{1}{c|}{$n=4$} \\
[.1cm]\hline
&&&&\\[-0.25cm]
1/3 & ${\bf E}(x)$ & $-1,81418881\times 10^{2}$ & $-1.08294258\times 10^3$ & $-3.08231679\times 10^3$ \\
& ${\bf H}(x)$ & $+0.34241316$ & $+0.17280892$ & $+0.34497729$\\ \cline{3-5}
&&&&\\[-0.3cm]         
& ${\bf E}(x)+{\bf H}(x)$ &$-1.81076468\times 10^{2}$  & $-1.08276977\times 10^3$ & $-3.08197181\times 10^3$\\
& $F_{n,\sa}(x;\mu)$ & $-1.80709370\times 10^{2}$  & $-1.08284759\times 10^{3}$ & $-3.08254767\times 10^{3}$ \\
\hline
&&&&\\[-0.25cm]
1/2 & ${\bf E}(x)$& $+0.06317153$ & $+1.15957937\times 10^3$ & $-4.47945373\times 10^4$ \\
& ${\bf H}(x)$&$+0.74012019$ &  $+1.09449277$ & $+1.45169481$ \\ \cline{3-5}
&&&&\\[-0.3cm]         
& ${\bf E}(x)+{\bf H}(x)$ &$+0.80329172$  & $+1.16067387\times 10^3$ & $-4.47930856\times 10^4$\\
& $F_{n,\sa}(x;\mu)$ & $+0.80329527$  & $+1.16069221\times 10^3$ & $-4.47921506\times 10^4$\\
\hline
&&&&\\[-0.25cm]
5/9 & ${\bf E}(x)$ &$\hspace{0.8cm} -\!\!-$ & $-0.14805870$ & $+2.77243091\times 10^2$\\
& ${\bf H}(x)$ &$+0.79825166$ & $+1.17615555$ & $+1.55857242$\\ \cline{3-5}
&&&&\\[-0.3cm]         
& ${\bf E}(x)+{\bf H}(x)$ &$+0.79825166$ & $+1.02809685$ & $+2.78801663\times 10^2$\\
& $F_{n,\sa}(x;\mu)$ &$+0.79825119$ & $+1.02809649$ & $+2.78801134\times 10^2$\\
\hline
&&&&\\[-0.25cm]
%2/3& ${\bf E}(x)$ & $\hspace{0.8cm} -\!\!-$ & $\hspace{0.8cm} -\!\!-$ & $\hspace{0.8cm} -\!\!-$ \\
2/3& ${\bf H}(x)$ & $+0.84046066$ & $+1.23266920$ & $+1.63072031$ \\ %\cline{3-5}
%&&&&\\[-0.3cm]         
%& ${\bf E}(x)+{\bf H}(x)$ &$+0.84046066$ & $+1.23266920$ & $+1.63072031$ \\
& $F_{n,\sa}(x;\mu)$ &$+0.84046066$  & $+1.23266920$ & $+1.63072031$\\
\hline
\end{tabular}
\end{center}
\end{table}

In Table 1 we present the results of numerical calculations for $x\to+\infty$ compared with the expansions given in Theorem 1. We choose four representative values of $\sa$ that focus on the different cases of Theorem 1 and $n=2, 3$ and 4.
The exact value of $F_{n,\sa}(x;\mu)$ was obtained by high-precision evaluation of (\ref{e41}). The exponential expansion ${\bf E}(x)$ was computed with truncation index $j=3$ and the algebraic expansion ${\bf H}(x)$ was optimally truncated (that is, at or near its smallest term). The first case $\sa=\f{1}{3}$ has an exponentially large expansion with a subdominant algebraic contribution for all three values of $n$. The second case $\sa=\fs$ corresponds to $n_0=2$; when $n=2$, ${\bf E}(x)$ is oscillatory and makes a similar contribution as ${\bf H}(x)$, whereas when $n=3$ and 4, ${\bf E}(x)$ is exponentially large. The third case $\sa=\f{5}{9}$ corresponds to $n_0=3$; when $n=2$ there is no exponential contribution, whereas when $n=3$, ${\bf E}(x)$ is oscillatory and so makes a similar contribution as ${\bf H}(x)$; when $n=4$, ${\bf E}(x)$ is exponentially large. Finally, when $\sa=\f{2}{3}$ the expansion of $F_{n,\sa}(x;\mu)$ is purely algebraic in character. 

\begin{table}[b]
\caption{\footnotesize{The values of the exponential and algebraic expansions compared with $F_{n,\sa}(x;\mu)$ for large $x<0$ for different values of $\sa$ and $n$ when $\mu=3/4$ and $|x|=8$ (for $\sa=1/4, 1/2, 2/5$), $|x|=5$ (for $\sa=3/4$).}}
\begin{center}
\begin{tabular}{|c r|l|l|l|}
\hline
&&&&\\[-0.3cm]
%\mcol{1}{|c|}{} & \mcol{1}{c|}{$\gamma=0$} &\mcol{1}{c|}{$\gamma=2$}& \mcol{1}{c|}{$\gamma=-2$} \\
\mcol{1}{|c}{$\sa$} & \mcol{1}{c|}{} & \mcol{1}{c|}{$n=2$} &\mcol{1}{c|}{$n=3$}& \mcol{1}{c|}{$n=4$} \\
[.1cm]\hline
&&&&\\[-0.25cm]
1/4 & ${\hat{\bf E}}(x)$ & $+1.59003829\times 10^{-2}$ & $+1.77442984\times 10^{-1}$ & $+6.49578248\times 10^{-1}$ \\
& $F_{n,\sa}(-x;\mu)$ & $+1.59003416\times 10^{-2}$  & $+1.77011100\times 10^{-1}$ & $+6.49580223\times 10^{-1}$ \\
\hline
&&&&\\[-0.25cm]
2/5 & ${\hat{\bf E}}(x)$ &$-4.18901636\times 10^{-2}$ & $-3.79446870\times 10^{0}$ & $-3.02428770\times 10^1$\\
& $F_{n,\sa}(-x;\mu)$ &$-4.18889220\times 10^{-2}$ & $-3.79475882\times 10^0$ & $-3.02402120\times 10^1$\\
\hline
&&&&\\[-0.25cm]
1/2 & ${\hat{\bf E}}(x)$& $-0.56022532$ & $+1.23070020\times 10^3$ & $-1.28808653\times 10^4$ \\
& $F_{n,\sa}(-x;\mu)$ & $-0.56023534$  & $+1.23066913\times 10^3$ & $-1.28803505\times 10^4$\\
\hline
&&&&\\[-0.25cm]
3/4& ${\hat{\bf E}}(x)$ & $+1.81213632\times 10^{28}$ & $+7.55354383\times 10^{13}$ & $-0.84956415$ \\
& ${\hat{\bf H}}(x)$ & $\hspace{0.8cm} -\!\!-$ & $-1.93112636\times 10^{-1}$ & $-0.28756658$ \\ \cline{3-5}
&&&&\\[-0.3cm]         
& ${\hat{\bf E}}(x)+{\hat{\bf H}}(x)$ &$+1.81213632\times 10^{28}$ & $+7.55354383\times 10^{13}$ & $-1.13713072$ \\
& $F_{n,\sa}(-x;\mu)$ &$+1.81213650\times 10^{28}$  & $+7.55354314\times 10^{13}$ & $-1.13713081$\\
\hline
\end{tabular}
\end{center}
\end{table}

In Table 2 we present illustrative examples of Theorem 2 when $x\to-\infty$. The first case $\sa=\f{1}{4}$ ($\ka=\f{3}{4}$),
has an expansion that is exponential in character; for $n<1/\sa=4$, ${\hat{\bf E}}(x)$ is exponentially small, whereas for $n=4$ the argument $\pi\ka-\om_0=\f{3}{8}\pi$ lies on the upper boundary of the exponentially large sector $|\arg\,z|<\f{3}{8}\pi$ and so ${\hat{\bf E}}(x)$ is oscillatory. For $n\geq 5$, ${\hat{\bf E}}(x)$ becomes exponentially large as $x\to-\infty$. In the second case $\sa=\f{2}{5}$ ($\ka=\f{3}{5}$), ${\hat{\bf E}}(x)$ is exponentially small for $n=2$ and exponentially large for $n\geq 3$. In the third case $\sa=\ka=\fs$, ${\hat{\bf E}}(x)$ is oscillatory for $n=2$ and exponentially large for $n\geq 3$. Finally, when $\sa=\f{3}{4}$ ($\ka=\f{1}{4}$)
the function $F_{n,\sa}(x;\mu)$ is exponentially large for $n=2, 3$ and $n\geq 5$. But for $n=4$, the two values 
$\om_0=\f{3}{8}\pi$, $\om_1=\f{1}{8}\pi$ yield arguments $\pi\ka-\om_r$ ($r=0, 1$) situated on {\it both\/} boundaries of the exponentially large sector $|\arg\,z|<\f{1}{8}\pi$. In this case ${\hat{\bf E}}(x)$ is oscillatory and, since $n^*=2$, there is in addition an algebraic contribution ${\hat{\bf H}}(x)$. 

It is seen from Tables 1 and 2 that the asymptotic values agree well with the numerically computed values of $F_{n,\sa}(\pm x;\mu)$, thereby confirming the accuracy of the expansions in Theorems 1 and 2.


\vspace{0.6cm}


\begin{thebibliography}{99}
%\footnotesize{
\bibitem{K}
L.L. Karasheva, On properties of an entire function that is a generalization of the Wright function, J. Math. Sciences {\bf 250} (2020) 753--759.

\bibitem{M}
F. Mainardi, {\it Fractional Calculus and Waves in Linear Viscoelasticity}, Imperial College Press, London, 2010.

\bibitem{DLMF}
F.W.J. Olver, D.W. Lozier, R.F. Boisvert and C.W. Clark (eds.),    
{\it NIST Handbook of Mathematical Functions}, Cambridge University Press, Cambridge, 2010.

\bibitem{P92}
R.B. Paris, Smoothing of the Stokes phenomenon using Mellin-Barnes integrals, J. Comput. Appl. Math. {\bf 41} (1992) 117--133. 

\bibitem{P1}
R.B. Paris, The asymptotics of the generalised Bessel function, Math. Aeterna {\bf 7} (2017) 381--406.

\bibitem{P2}
R,B, Paris, Asymptotics of the special functions of fractional calculus, in {\it Handbook of Fractional Calculus with Applications}, vol. 1 (eds. A. Kochubei and Y. Luchko) pp.~297--325, De Gruyter, Berlin, 2019.

\bibitem{PK} 
R.B. Paris and D. Kaminski,  {\it Asymptotics and Mellin-Barnes Integrals}, 
Cambridge University Press, Cambridge, 2001.

%\bibitem{WW} 
%E. T. Whittaker and G. N. Watson, {\it Modern Analysis}, Cambridge University Press, Cambridge, 1952.

%}
\end{thebibliography}
\end{document}